# DECONVOLUTION WITH UNKNOWN ERROR DISTRIBUTION

By Jan Johannes

*Ruprecht–Karls–Universität Heidelberg*

We consider the problem of estimating a density $f_X$ using a sample $Y_1, \ldots, Y_n$ from $f_Y = f_X \star f_\epsilon$, where $f_\epsilon$ is an unknown density. We assume that an additional sample $\epsilon_1, \ldots, \epsilon_m$ from $f_\epsilon$ is observed. Estimators of $f_X$ and its derivatives are constructed by using nonparametric estimators of $f_Y$ and $f_\epsilon$ and by applying a spectral cut-off in the Fourier domain. We derive the rate of convergence of the estimators in case of a known and unknown error density $f_\epsilon$, where it is assumed that $f_X$ satisfies a polynomial, logarithmic or general source condition. It is shown that the proposed estimators are asymptotically optimal in a minimax sense in the models with known or unknown error density, if the density $f_X$ belongs to a Sobolev space $H_\mathbb{p}$ and $f_\epsilon$ is ordinary smooth or supersmooth.

**1. Introduction.** Let $X$ and $\epsilon$ be independent random variables with unknown density functions $f_X$ and $f_\epsilon$, respectively. The objective is to nonparametrically estimate the density function $f_X$ and its derivatives based on a sample of $Y = X + \epsilon$. In this setting, the density $f_Y$ of $Y$ is the convolution of the density of interest, $f_X$, and the density $f_\epsilon$ of the additive noise, that is,

$$(1.1) \qquad f_Y(y) = f_X \star f_\epsilon(y) := \int_{-\infty}^{\infty} f_X(x) f_\epsilon(y - x) \, dx.$$

Suppose we observe $Y_1, \ldots, Y_n$ from $f_Y$ and the error density $f_\epsilon$ is known. Then, the estimation of the deconvolution density $f_X$ is a classical problem in statistics. The most popular approach is to estimate $f_Y$ by a kernel estimator and then solve (1.1) using a Fourier transform (see Carroll and Hall [4], Devroye [7], Efromovich [9], Fan [11, 12], Stefanski [36], Zhang [41], Goldenshluger [[14], [15]] and Kim and Koo [21]). Spline-based methods are considered, for example, in Mendelsohn and Rice [28] and Koo and Park [22].

The estimation of the deconvolution density using a wavelet decomposition is studied in Pensky and Vidakovic [34], Fan and Koo [13] and Bigot and Van Bellegem [1], while Hall and Qiu [16] have proposed a discrete Fourier series expansion. A penalization and projection approach is proposed in Carrasco and Florens [3] and Comte, Rozenholc and Taupin [6].

The underlying idea behind all approaches is to replace in (1.1) the unknown density $f_Y$ by its estimator and then solve (1.1). However, solving (1.1) leads to an ill-posed inverse problem and, hence, the inversion of (1.1) has to be "regularized" in some way. We now describe three examples of regularization. The first example is kernel estimators, where the kernel has a limited bandwidth, that is, the Fourier transform of the kernel has a bounded support. In this case, asymptotic optimality, both pointwise and global, over a class of functions whose derivatives are Lipschitz continuous, is proven in Carroll and Hall [4] and Fan [11, 12]. The second example is estimators based on a wavelet decomposition, where the wavelets have limited bandwidths. For the wavelet estimator, Pensky and Vidakovic [34] show asymptotic optimality of the mean integrated squared error (MISE) over the Sobolev space $H_p$, which describes the level of smoothness of a function $f$ in terms of its Fourier transform $\mathcal{F}f$. In the third example, the risk in the Sobolev norm of $H_s$ ($H_s$-risk) and asymptotic optimality over $H_p$, $p \geq s$, of an estimator using a spectral cut-off (thresholding of the Fourier transform $\mathcal{F}f_\epsilon$ of $f_\epsilon$) is derived in Mair and Ruymgaart [26].

However, in the above examples, $f_X$ and $f_\epsilon$ are assumed to be ordinary smooth or supersmooth, that is, their Fourier transforms have polynomial or exponential descent. All these cases can be characterized by a "source condition" (defined below), which allows for more general tail behavior.

In several applications, for example, in optics and medicine (cf. Tessier [38] and Levitt [23]), the noise density $f_\epsilon$ may be unknown. In this case, without any additional information, the density $f_X$ cannot be recovered from the density of $f_Y$ through (1.1), that is, the density $f_X$ is not identified if only a sample $Y_1, \ldots, Y_n$ from $f_Y$ is observed. It is worth noting that in some special cases the deconvolution density $f_X$ can be identified (cf. Butucea and Matias [2] or Meister [27]). Deconvolution without prior knowledge of the error distribution is also possible in the case of panel data (cf. Horowitz and Markatou [19], Hall and Yao [17] or Neumann [32]).

In this paper, we deal with the estimation of a deconvolution density $f_X$ when only an approximation of the error density $f_\epsilon$ is given. More precisely, following Diggle and Hall [8] we suppose, that in addition to a sample $Y_1, \ldots, Y_n$ from $f_Y$, we observe a sample $\epsilon_1, \ldots, \epsilon_m$ from $f_\epsilon$. An interesting example in bio-informatics can be found in the analysis of cDNA microarrays, where $Y$ is the intensity measure, $X$ is the expressed gene intensity and $\epsilon$ is the background intensity (for details see Havilio [18]). In a situation where an estimator of $f_\epsilon$ is used, rather than the true density, Neumann [31] shows asymptotic optimality of the MISE over the Bessel-potential space when the error density is ordinary smooth. In case of a circular convolution problem, Cavalier and Hengartner [5] present oracle inequalities and adaptive estimation. However, they also assume the error density to be ordinary smooth. By constraining the error density to be ordinary smooth, a rich class of distributions, such as the normal distribution, are excluded. The purpose of this paper is to propose and study a deconvolution scheme which has enough flexibility to allow a wide range of tail behaviors of $\mathcal{F}f_X$ and $\mathcal{F}f_\varepsilon$.

The estimators of the deconvolution density considered in this paper are based on a regularized inversion of (1.1) using a spectral cut-off, where we replace the unknown density $f_Y$ by a nonparametric estimator and the Fourier transform of $f_\epsilon$ by its empirical counterpart. We derive the $H_s$-risk of the proposed estimator for a wide class of density functions, which unifies and generalizes many of the previous results for known and unknown error density. Roughly speaking, we show in case of known $f_\epsilon$ that the $H_s$-risk can be decomposed into a function of the MISE of the nonparametric estimator of $f_Y$ plus an additional bias term which is a function of the threshold (the parameter which determines the spectral cut-off point). The relationship between $\mathcal{F}f_X$ and $\mathcal{F}f_\epsilon$ is then essentially determining the functional form of the bias term. For example, the bias is a logarithm of the threshold when the error distribution is supersmooth (e.g., normal) and $f_X$ is ordinary smooth (e.g., double exponential). On the other hand, if both the error distribution and $f_X$ are ordinary smooth or supersmooth, the bias is a polynomial of the threshold. We show that the theory behind these rates can be unified using an index function $\kappa$ (cf. Nair, Pereverzev and Tautenhahn [29]), which "links" the tail behavior of $\mathcal{F}f_X$ and $\mathcal{F}f_\epsilon$ by supposing that $|\mathcal{F}f_X|^2/\kappa(|\mathcal{F}f_\varepsilon|^2)$ is integrable.

Under certain conditions on the index function, we prove that the $H_s$-risk in the model with unknown $f_\epsilon$ can be decomposed into a part with the same bound as the $H_s$-risk for known $f_\epsilon$ and a second term which is only a function of the sample size $m$ (of errors $\epsilon$). The functional form of the second term is then again determined by the relationship between $\mathcal{F}f_X$ and $\mathcal{F}f_\epsilon$. We show that the second term provides a lower bound for the $H_s$-risk on its own and, hence, cannot be avoided. It follows that the estimator is minimax in the model with unknown $f_\epsilon$ when the bound of the $H_s$-risk for known $f_\epsilon$ is of minimax optimal order. Furthermore, it is of interest to compare the rates of convergence of the $H_s$-risk when the density of $f_\epsilon$ is estimated with the rates, where $f_\epsilon$ is known. We show that under certain conditions on the index function, a sample size $m$ which increases at least as fast as the inverse of the MISE of the nonparametric estimator of $f_Y$, ensures an asymptotically negligible estimation error of $f_\epsilon$. However, in special cases even slower rates of $m$ are enough.

In this paper, we use the classical Rosenblatt–Parzen kernel estimator (cf. Parzen [33]) without a limited bandwidth to estimate the density $f_Y$. However, since the $H_s$-risk of the proposed estimator can be decomposed using the MISE of the density estimator of $f_Y$, any other nonparametric estimation method (e.g., based on splines or wavelets) can be used and the theory still holds.

The paper is organized in the following way. In Section 2, we give a brief description of the background of the methodology and we define the estimator of $f_X$ when the density $f_\epsilon$ is known as well as when $f_\epsilon$ is unknown. We investigate the asymptotic behavior of the estimator of $f_X$ in case of a known and an unknown density $f_\epsilon$ in Sections 3 and 4, respectively. All proofs can be found in the Appendix.

## 2. Methodology.

2.1. *Background to methodology.* In this paper, we suppose that $f_X$ and $f_\epsilon$ [hence also $f_Y$] are contained in the set $\mathcal{D}$ of all densities in $L^2(\mathbb{R})$, which is endowed with the usual norm $\|\cdot\|$. We use the notation $[\mathcal{F}g](t)$ for the Fourier transform $\frac{1}{\sqrt{2\pi}} \int_{-\infty}^{\infty} \exp(-itx) g(x)\,dx$ of a function $g \in L^1(\mathbb{R}) \cap L^2(\mathbb{R})$, which is unitary. Since $X$ and $\epsilon$ are assumed to be independent, the Fourier transform of $f_Y$ satisfies $\mathcal{F}f_Y = \sqrt{2\pi} \cdot \mathcal{F}f_X \cdot \mathcal{F}f_\epsilon$. Therefore, assuming $|[\mathcal{F}f_\epsilon](t)|^2 > 0$, for all $t \in \mathbb{R}$, the density $f_X$ can be recovered from $f_Y$ and $f_\epsilon$ by

$$(2.1) \qquad \mathcal{F}f_X = \frac{\mathcal{F}f_Y \cdot \overline{\mathcal{F}f_\epsilon}}{\sqrt{2\pi} \cdot |\mathcal{F}f_\epsilon|^2},$$

where $\overline{\mathcal{F}f_\epsilon}$ denotes the complex conjugate of $\mathcal{F}f_\epsilon$. It is well known that replacing in (2.1) the unknown density $f_Y$ by a consistent estimator $\widehat{f_Y}$ does not in general lead to a consistent estimator of $f_X$. To be more precise, since $|\mathcal{F}f_\epsilon|^{-1}$ is not bounded, $\mathbb{E}\|\widehat{f_Y} - f_Y\|^2 = o(1)$ does not generally imply $\mathbb{E}\|[\mathcal{F}\widehat{f_Y} - \mathcal{F}f_Y] \cdot |\mathcal{F}f_\epsilon|^{-1}\|^2 = o(1)$, that is, the inverse operation of a convolution is not continuous. Therefore, the deconvolution problem is ill posed in the sense of Hadamard. In the literature, several approaches are proposed in order to circumvent this instability issue. Essentially, all of them replace (2.1) with a regularized version that avoids having the denominator becoming too small [e.g., nonparametric methods using a kernel with limited bandwidth estimate $\mathcal{F}f_Y(t)$, and also $\mathcal{F}f_X(t)$, for $|t|$ larger than a threshold by zero]. There are a large number of alternative regularization schemes in the numerical analysis literature available, such as the Tikhonov regularization, Landweber iteration or the $\nu$-methods, to name but a few (cf. Engl, Hanke and Neubauer [10]). However, in this paper we regularize (2.1) by introducing a threshold $\alpha > 0$ and a function $\ell_s(t) := (1 + t^2)^{s/2}$,

$s, t \in \mathbb{R}$, that is, for $s \geq 0$, we consider the regularized version $f_{X_s}^\alpha$ given by

$$(2.2) \qquad \mathcal{F} f_{X_s}^\alpha := \frac{\mathcal{F} f_Y \cdot \overline{\mathcal{F} f_\epsilon}}{\sqrt{2\pi} \cdot |\mathcal{F} f_\epsilon|^2} \cdot \mathbb{1}\{|\mathcal{F} f_\epsilon / \ell_s|^2 \geq \alpha\}.$$

Then, $f_{X_s}^\alpha$ belongs to the well-known Sobolev space $H_s$ defined by

$$(2.3) \qquad H_s := \left\{ f \in L^2(\mathbb{R}) : \|f\|_s^2 := \int_{-\infty}^\infty (1+t^2)^s |[\mathcal{F} f](t)|^2 \, dt < \infty \right\}.$$

Moreover, let $H_s^\rho := \{f \in H_s : \|f\|_s^2 \leq \rho\}$, for $\rho > 0$. Thresholding in the Fourier domain has been used, for example, in Devroye [7], Liu and Taylor [24], Mair and Ruymgaart [26] or Neumann [31] and coincides with an approach called spectral cut-off in the numerical analysis literature (cf. Tautenhahn [37]).

2.2. *Estimation of $f_X$ when $f_\epsilon$ is known.* Let $Y_1, \ldots, Y_n$ be an i.i.d. sample of $Y$, which we use to construct an estimator $\widehat{f_Y}$ of $f_Y$. The estimator $\widetilde{f_{X_s}}$ of $f_X$ based on the regularized version (2.2) is then defined by

$$(2.4) \qquad \mathcal{F} \widetilde{f_{X_s}} := \frac{\mathcal{F} \widehat{f_Y} \cdot \overline{\mathcal{F} f_\epsilon}}{\sqrt{2\pi} \cdot |\mathcal{F} f_\epsilon|^2} \cdot \mathbb{1}\{|\mathcal{F} f_\epsilon / \ell_s|^2 \geq \alpha\},$$

where the threshold $\alpha := \alpha(n)$ has to tend to zero as the sample size $n$ increases. The truncation in the Fourier domain will lead as usual to a bias term which is a function of the threshold. In Lemma A.1 in the Appendix, we show that by using this specific structure for the truncation, the functional form of the bias term is determined by the relationship between $\mathcal{F} f_X$ and $\mathcal{F} f_\epsilon$. In this paper, we stick to a nonparametric kernel estimation approach, but we would like to stress that any other density estimation procedure could be used as well. The kernel estimator of $f_Y$ is defined by

$$(2.5) \qquad \widehat{f_Y}(y) := \frac{1}{nh} \sum_{j=1}^n K\left(\frac{Y_j - y}{h}\right), \qquad y \in \mathbb{R},$$

where $h > 0$ is a bandwidth and $K$ a kernel function. As usual in the context of nonparametric kernel estimation the bandwidth $h$ has to tend to zero as the sample size $n$ increases. In order to derive a rate of convergence of $\widehat{f_Y}$, we follow Parzen [33] and consider, for each $r \geq 0$, the class of kernel functions

$$(2.6) \qquad \mathcal{K}_r := \left\{ K \in L^1(\mathbb{R}) \cap L^2(\mathbb{R}) : \lim_{t \to 0} \frac{|1 - \sqrt{2\pi}[\mathcal{F} K](t)|}{|t|^r} = \kappa_r < \infty \right\}.$$

If $f_Y \in H_r^q$, for $q, r > 0$, then the MISE of the estimator $\widehat{f_Y}$ given in (2.5), constructed by using a kernel $K \in \mathcal{K}_r$ and a bandwidth $h = cn^{-1/(2r+1)}$, $c > 0$, is of order $n^{-2r/(2r+1)}$ (cf. Parzen [33]) and, hence, obtains the minimax optimal order over the class $H_r^q$ (cf. [40], Chapter 24).

2.3. *Estimation of $f_X$ given an estimator of $f_\epsilon$.* Suppose $Y_1, \ldots, Y_n$ and $\epsilon_1, \ldots, \epsilon_m$ form i.i.d. samples of $f_Y$ and $f_\epsilon$, respectively. We consider again the nonparametric kernel estimator $\widehat{f_Y}$ defined in (2.5). In addition, we estimate the Fourier transform $\mathcal{F}f_\epsilon$ using its empirical counterpart, that is,

$$(2.7) \qquad [\widehat{\mathcal{F}f_\epsilon}](t) := \frac{1}{m \cdot \sqrt{2\pi}} \sum_{j=1}^{m} e^{-it\epsilon_j}, \qquad t \in \mathbb{R}.$$

Then, the estimator $\widehat{f_X}_s$ based on the regularized version (2.2) is defined by

$$(2.8) \qquad \mathcal{F}\widehat{f_X}_s := \frac{\mathcal{F}\widehat{f_Y} \cdot \overline{\widehat{\mathcal{F}f_\epsilon}}}{\sqrt{2\pi} \cdot |\widehat{\mathcal{F}f_\epsilon}|^2} \cdot \mathbb{1}\{|\widehat{\mathcal{F}f_\epsilon}/\ell_s|^2 \geq \alpha\},$$

where $\alpha := \alpha(n, m)$ has to tend to zero as the sample sizes $n$ and $m$ increase.

**3. Theoretical properties of the estimator when $f_\epsilon$ is known.** We shall measure the performance of the estimator $\widetilde{f_X}_s$ defined in (2.4) by the $H_s$-risk, that is, $\mathbb{E}\|\widetilde{f_X}_s - f_X\|_s^2$, provided $f_X \in H_p$, for some $p \geq s \geq 0$. For an integer $k$, the Sobolev norm $\|g\|_k$ is equivalent to $\|g\| + \|g^{(k)}\|$, where the $k$th weak derivative $g^{(k)}$ of $g$ satisfies $[\mathcal{F}g^{(k)}](t) := (-it)^k[\mathcal{F}g](t)$. Therefore, the $H_k$-risk reflects the performance of $\widetilde{f_X}_k$ and $\widetilde{f_X}_k^{(k)}$ as estimators of $f_X$ and $f_X^{(k)}$, respectively. However, in what follows a situation without an a priori assumption on the smoothness of $f_X$ is also covered considering $p = s = 0$.

The $H_s$-risk is essentially determined by the MISE of the estimator of $f_Y$ and by the regularization bias. To be more precise, by using $f_{Xs}^\alpha$ given in (2.2) and assuming $f_X \in H_p$, for some $p \geq s \geq 0$, we bound the $H_s$-risk by

$$(3.1) \qquad \mathbb{E}\|\widetilde{f_X}_s - f_X\|_s^2 \leq \pi^{-1}\alpha^{-1}\mathbb{E}\|\widehat{f_Y} - f_Y\|^2 + 2\|f_{Xs}^\alpha - f_X\|_s^2,$$

where, due to Lebesgue's dominated convergence theorem, the regularization bias satisfies $\|f_{Xs}^\alpha - f_X\|_s^2 = o(1)$ as $\alpha$ tends to zero.

PROPOSITION 3.1. *Suppose that $f_X \in H_p$, $p \geq 0$. Let $\widehat{f_Y}$ be a consistent estimator of $f_Y$, that is, $\mathbb{E}\|\widehat{f_Y} - f_Y\|^2 = o(1)$ as $n \to \infty$. Consider, for $0 \leq s \leq p$, the estimator $\widetilde{f_X}_s$ given in (2.4) with threshold satisfying $\alpha = o(1)$ and $\mathbb{E}\|\widehat{f_Y} - f_Y\|^2/\alpha = o(1)$ as $n \to \infty$. Then, $\mathbb{E}\|\widetilde{f_X}_s - f_X\|_s^2 = o(1)$ as $n \to \infty$.*

In order to obtain a rate of convergence of the regularization bias and, hence, the $H_s$-risk of $\widetilde{f_X}_s$, we consider first a *polynomial source condition*

$$(3.2) \qquad \rho := \|\ell_s \cdot \mathcal{F}f_X \cdot (|\mathcal{F}f_\epsilon/\ell_s|^2)^{-\beta/2}\| < \infty \qquad \text{for some } \beta > 0, s \geq 0.$$

Note that (3.2) implies that $f_X \in H_s$.

EXAMPLE 3.1. To illustrate this and also the following source conditions, let us consider three different types of densities. These are, (i) the density $g$ of a symmetrized $\chi^2$ distribution with $k$ degrees of freedom, that is, $[\mathcal{F}g](t) = (2\pi)^{-1/2}(1+4t^2)^{-k/2}$, (ii) the density $g$ of a centered Cauchy distribution with scale parameter $\gamma > 0$, that is, $[\mathcal{F}g](t) = (2\pi)^{-1/2}\exp(-\gamma|t|)$, and (iii) the density $g$ of a centered normal distribution with variance $\sigma^2 > 0$, that is, $[\mathcal{F}g](t) = (2\pi)^{-1/2}\exp(-\sigma^2 t^2/2)$. Suppose $f_X$ and $f_\epsilon$ are symmetrized $\chi^2$ densities with $k_X$ and $k_\epsilon$ degrees of freedom, respectively. Then, the polynomial source condition (3.2) is only satisfied for $0 \leq s < k_X - 1/2$. If $f_X$ and $f_\epsilon$ are Cauchy densities or $f_X$ and $f_\epsilon$ are Gaussian densities, then $\mathcal{F}f_X$ and $\mathcal{F}f_\epsilon$ descend exponentially and (3.2) holds for all $s \geq 0$.

THEOREM 3.2. *Suppose that $f_X$ satisfies the polynomial source condition (3.2), for some $s \geq 0$ and $\beta > 0$. Consider the estimator $\widetilde{f_X}_s$ defined in (2.4) by using a threshold $\alpha = c \cdot (\mathbb{E}\|\widehat{f_Y} - f_Y\|^2)^{1/(\beta+1)}$, $c > 0$. Then, there exists a constant $C > 0$ depending only on $\rho$ given in (3.2), $\beta$ and $c$ such that $\mathbb{E}\|\widetilde{f_X}_s - f_X\|_s^2 \leq C \cdot (\mathbb{E}\|\widehat{f_Y} - f_Y\|^2)^{\beta/(\beta+1)}$, as $\mathbb{E}\|\widehat{f_Y} - f_Y\|^2 \to 0$.*

REMARK 3.1. In Lemma A.1 in the Appendix, we show by applying standard techniques for regularization methods that the polynomial source condition (3.2) implies $\|f_{X_s}^\alpha - f_X\|_s^2 \leq \alpha^\beta \rho^2$. Then, we obtain the result by balancing in (3.1) the two terms on the right-hand side. On the other hand, from Theorem 4.11 in Engl, Hanke and Neubauer [10] follows that $\|f_{X_s}^\alpha - f_X\|_s^2 = O(\alpha^\eta)$, for some $\eta > 0$, implies (3.2) for all $\beta < \eta$, that is, the order $O(\alpha^\beta)$ is optimal over the class $\{f_X$ satisfies $(3.2)\}$. Therefore, one would expect that an optimal estimation of $f_Y$ leads to an optimal estimation of $f_X$. However, the polynomial source condition is not sufficient to derive an optimal rate of convergence of the MISE of $\widehat{f_Y}$ over the class $\{f_Y = f_\epsilon \star f_X : f_X$ satisfies $(3.2)\}$. For example, if $f_\epsilon$ is a Gaussian density, this class contains only analytic functions, while it equals $H_{(\beta+1)(s+1)}$ when $f_\epsilon$ is a Laplace density.

Without further information about $f_\epsilon$ it is difficult to give for arbitrary $\beta > 0$ an interpretation of the polynomial source condition. However, if we suppose additionally that $f_\epsilon$ is *ordinary smooth*, that is, there exists $a > 1/2$ and a constant $d > 0$, such that

(3.3) $\qquad d \leq (1+t^2)^a |[\mathcal{F}f_\epsilon](t)|^2 \leq d^{-1} \qquad$ for all $t \in \mathbb{R}$.

Then, the smoothness condition $f_X \in H_p$, for some $p > 0$, is equivalent to the polynomial source condition (3.2) with $0 \leq s < p$ and $\beta = (p-s)/(s+a)$. Moreover, we have $H_{p+a} = \{f_Y = f_\epsilon \star f_X : f_X \in H_p\}$, for all $p \geq 0$. Therefore, the convolution with $f_\epsilon$ is also called *finitely smoothing* (cf. Mair and Ruymgaart [26]). From Theorem 3.2, we obtain the following corollary, which establishes the optimal rate of convergence of $\widetilde{f_X}_s$ over $H_p$.

COROLLARY 3.3. *Suppose that $f_X \in H_p$, $p > 0$ and $f_\epsilon$ satisfies (3.3) for $a > 1/2$. Let $\widehat{f_Y}$ defined in (2.5) be constructed using a kernel $K \in \mathcal{K}_{p+a}$ [see (2.6)] and a bandwidth $h = cn^{-1/(2(p+a)+1)}$, $c > 0$. Consider for $0 \leq s < p$ the estimator $\widetilde{f_X}_s$ defined in (2.4) with threshold $\alpha = cn^{-2(a+s)/(2(p+a)+1)}$, $c > 0$. Then, we have $\mathbb{E}\|\widetilde{f_X}_s - f_X\|_s^2 = O(n^{-2(p-s)/(2(a+p)+1)})$ as $n \to \infty$.*

REMARK 3.2. The rate of convergence in the last result is known to be minimax optimal over the class $H_p^\rho$, provided that the density $f_\epsilon$ satisfies (3.3) (cf. Mair and Ruymgaart [26]). Since under the assumptions of the corollary $f_X$ belongs to $H_p$ if and only if $f_Y$ lies in $H_{p+a}$, it follows that the kernel estimator of $f_Y$ is constructed such that its MISE has the minimax optimal order over the class $H_{p+a}^q$. Moreover, using an estimator of $f_Y$ which does not have an order optimal MISE, the estimator of $f_X$ would not reach the minimax optimal rate of convergence. Hence, in this situation the optimal estimation of $f_Y$ is necessary to obtain an optimal estimator of $f_X$. We shall emphasize the role of the parameter $a$, which specifies through the condition (3.3) the tail behavior of the Fourier transform $\mathcal{F}f_\epsilon$. As we see, if the value $a$ increases, the obtainable optimal rate of convergence decreases. Therefore, the parameter $a$ is often called *degree of ill posedness* (cf. Natterer [30]).

If, for example, $f_X$ is a Laplace and $f_\epsilon$ is a Cauchy or Gaussian density, then not a polynomial but a *logarithmic source condition* holds true, that is,

$$(3.4) \quad \rho := \|\ell_s \cdot \mathcal{F}f_X \cdot |\ln(|\mathcal{F}f_\epsilon/\ell_s|^2)|^{\beta/2}\| < \infty \qquad \text{for some } \beta > 0,\ s \geq 0.$$

THEOREM 3.4. *Let $f_X$ satisfy the logarithmic source condition (3.4), for some $s \geq 0$ and $\beta > 0$. Consider the estimator $\widetilde{f_X}_s$ defined in (2.4) by using a threshold $\alpha = c \cdot (\mathbb{E}\|\widehat{f_Y} - f_Y\|^2)^{1/2}$, for some $c > 0$. Then, there exists a constant $C > 0$ depending only on $\rho$ given in (3.4), $\beta$ and $c$ such that we have $\mathbb{E}\|\widetilde{f_X}_s - f_X\|_s^2 \leq C \cdot |\log(\mathbb{E}\|\widehat{f_Y} - f_Y\|^2)|^{-\beta}$, as $\mathbb{E}\|\widehat{f_Y} - f_Y\|^2 \to 0$.*

Additionally, if we assume that the density $f_\epsilon$ is *supersmooth*, that is, there exists $a > 0$ and a constant $d > 0$, such that

$$(3.5) \qquad d \leq (1+t^2)^a |\ln(|[\mathcal{F}f_\epsilon](t)|^2)|^{-1} \leq d^{-1} \qquad \text{for all } t \in \mathbb{R},$$

then the smoothness condition $f_X \in H_p$, $p > 0$ is equivalent to the logarithmic source condition (3.4), with $0 \leq s < p$ and $\beta = (p-s)/a$. Moreover, $f_\epsilon$, and therefore $f_Y$, belong to $H_r$, for all $r > 0$, and given $a \geq 1$, $f_\epsilon$ and hence $f_Y$, are analytic functions (cf. Kawata [20]). Therefore, the convolution with $f_\epsilon$ is called *infinitely smoothing* (cf. Mair and Ruymgaart [26]).

COROLLARY 3.5. *Suppose that $f_X \in H_p$, $p > 0$ and $f_\epsilon$ satisfies (3.5) for some $a > 0$. Let $\widehat{f_Y}$ given in (2.5) be constructed by using a kernel $K \in \mathcal{K}_r$ [see (2.6)] and a bandwidth $h = cn^{-1/(2r+1)}$, $c, r > 0$. Consider, for $0 \leq s < p$, the estimator $\widetilde{f_X}_s$ defined in (2.4) with threshold $\alpha = cn^{-r/(2r+1)}$, $c > 0$. Then, we have $\mathbb{E}\|\widetilde{f_X}_s - f_X\|_s^2 = O((\log n)^{-(p-s)/a})$, as $n \to \infty$.*

REMARK 3.3. The rate of convergence in Corollary 3.5 is again minimax optimal over the class $H_p^\rho$, given that the density $f_\epsilon$ satisfies (3.5) (cf. Mair and Ruymgaart [26]). It seems rather surprising that in opposite to Corollary 3.3, an increasing value $r$ improves the order of the MISE of the estimator $\widehat{f_Y}$ uniform over the class $\{f_Y = f_\epsilon \star f_X : f_X \in H_p^\rho\}$, but does not change the order of the $H_s$-risk of $\widetilde{f_X}_s$ (compare Remark 3.2). This, however, is due to the fact that the $H_s$-risk of $\widetilde{f_X}_s$ is of order $O(n^{-r/(2r+1)}) + O((\log n^{r/(2r+1)})^{-(p-s)/a}) = O((\log n)^{-(p-s)/a})$. So $r$ does not appear formally, but is actually hidden in the order symbol. Note that neither the bandwidth $h$ nor the threshold $\alpha$ depends on the level $p$ of smoothness of $f_X$, that is, the estimator is adaptive. Moreover, the parameter $a$ specifying in condition (3.5) the tail behavior of the Fourier transform $\mathcal{F}f_\epsilon$, in this situation also describes the *degree of ill posedness*.

Consider, for example, a Cauchy density $f_X$ and a Gaussian density $f_\epsilon$, then neither the polynomial source condition (3.2) nor the logarithmic source condition (3.4) is appropriate. However, both source conditions can be unified and extended using an index function $\kappa : (0, 1] \to \mathbb{R}^+$, which we always assume here to be a continuous and strictly increasing function with $\kappa(0+) = 0$ (cf. Nair, Pereverzev and Tautenhahn [29]). Then, we consider a *general source condition*

$$(3.6) \quad \rho := \|\ell_s \cdot \mathcal{F}f_X \cdot |\kappa(|\mathcal{F}f_\epsilon/\ell_s|^2)|^{-1/2}\| < \infty \qquad \text{for some } s \geq 0.$$

THEOREM 3.6. *Let $f_X$ satisfy the general source condition (3.6) for some concave index function $\kappa$ and $s \geq 0$. Denote by $\Phi$ and $\omega$ the inverse function of $\kappa$ and $\omega^{-1}(t) := t\Phi(t)$, respectively. Consider the estimator $\widetilde{f_X}_s$ defined in (2.4) by using $\alpha = c \cdot \mathbb{E}\|\widehat{f_Y} - f_Y\|^2 / \omega(c \cdot \mathbb{E}\|\widehat{f_Y} - f_Y\|^2)$, $c > 0$. Then, there exists a constant $C > 0$ depending only on $\rho$ given in (3.6) and $c$ such that $\mathbb{E}\|\widetilde{f_X}_s - f_X\|_s^2 \leq C \cdot \omega(\mathbb{E}\|\widehat{f_Y} - f_Y\|^2)$, as $\mathbb{E}\|\widehat{f_Y} - f_Y\|^2 \to 0$.*

REMARK 3.4. (i) Let $\mathcal{S}_{f_\epsilon}^\gamma$ be the set of all densities $f_X$ satisfying the general source condition (3.4) with $\rho \leq \gamma$. We define the modulus of continuity $\omega(\delta, \mathcal{S}_{f_\epsilon}^\gamma) := \sup\{\|g\|_s^2 : g \in \mathcal{S}_{f_\epsilon}^\gamma, \|f_\epsilon \star g\|^2 \leq \delta\}$ of the inverse operation of a convolution with $f_\epsilon$ over the set $\mathcal{S}_{f_\epsilon}^\gamma \subset H_s$. Since the index function $\kappa$ is assumed to be concave, it follows that the inverse function of $\omega$ is convex. Then, by using Theorem 2.2 in Nair, Pereverzev and Tautenhahn [29], we have $\omega(\delta) = O(\omega(\delta, \mathcal{S}_{f_\epsilon}^\gamma))$, as $\delta \to 0$. In the case of a deterministic approximation $f_Y^\delta$ of $f_Y$ with $\|f_Y^\delta - f_Y\| \leq \delta$, it is shown in Vainikko and Veretennikov [39] that $\omega(\delta, \mathcal{S}_{f_\epsilon}^\gamma)$ provides a lower bound over the class $\mathcal{S}_{f_\epsilon}^\gamma$ of the approximation error for any deconvolution method based only on $f_Y^\delta$. Therefore, we conjecture, that the bound in Theorem 3.6 is order optimal over the class $\mathcal{S}_{f_\epsilon}^\gamma$, given the MISE of $f_Y$ is order optimal over the class $\{f_Y = f_X \star f_\epsilon, f_X \in \mathcal{S}_{f_\epsilon}^\gamma\}$.

(ii) Define $\kappa(t) := |\log(ct)|^{-\beta}$, $c := \exp(-1 - \beta)$. Then, $\kappa$ is a concave index function and $\omega(\delta) = |\log \delta|^{-\beta}(1 + o(1))$, as $\delta \to 0$ (see Mair [25]). Thus, the result under a logarithmic source condition (Theorem 3.4) is covered by Theorem 3.6. However, the index function $\kappa(t) = t^\beta$ is concave only if $\beta \leq 1$, and hence the result in the case of a polynomial source condition (Theorem 3.2) is only partially obtained by Theorem 3.6. Nevertheless, we can apply Theorem 3.6 in the situation of a Cauchy density $f_X$ and a Gaussian density $f_\epsilon$ (compare Example 3.1), since in this case, for all $0 < \beta < 2\gamma/\sigma$ and $s \geq 0$, the general source condition is satisfied with concave index function $\kappa(t) = \exp(-\beta\sqrt{|\log(ct)|})$, $c := \exp(-(\beta^2 \vee 2))$. Moreover, if we denote $h(t) := (t/\beta + \beta/2)^2$, then $\omega^{-1}(t) = \exp(-h(-\log t))/c'$, with $c' = \exp(\beta^2/4 + (\beta^2 \vee 2))$. Since $\omega(t) = \exp(-h^{-1}(-\log t/c'))$, with $h^{-1}(y) = \beta\sqrt{y} - \beta^2/2$ for all $y \geq \beta^2/4$, we conclude that the $H_s$-risk in this case is of order $\exp(-\beta|\log \mathbb{E}\|\widehat{f_Y} - f_Y\|^2|^{1/2})$.

**4. Theoretical properties of the estimator when $f_\epsilon$ is unknown.** Let $\widehat{f_{Xs}^\alpha}$ be defined by $\mathcal{F}\widehat{f_{Xs}^\alpha} := \mathbb{1}\{|\widehat{\mathcal{F}f_\epsilon}/\ell_s|^2 \geq \alpha\} \cdot \mathcal{F}f_X$. Then, assuming $f_X \in H_p$, $p \geq s$, we bound the $H_s$-risk of $\widehat{f_{Xs}}$ given in (2.8) by

$$(4.1) \qquad \mathbb{E}\|\widehat{f_{Xs}} - f_X\|_s^2 \leq 2\mathbb{E}\|\widehat{f_{Xs}} - \widehat{f_{Xs}^\alpha}\|_s^2 + 2\mathbb{E}\|\widehat{f_{Xs}^\alpha} - f_X\|_s^2,$$

where we show in the proof of the next proposition that $\mathbb{E}\|\widehat{f_{Xs}} - \widehat{f_{Xs}^\alpha}\|_s^2$ is bounded up to a constant by $\alpha^{-1}(\mathbb{E}\|\widehat{f_Y} - f_Y\|^2 + m^{-1})$, and that the "regularization error" satisfies $\mathbb{E}\|\widehat{f_{Xs}^\alpha} - f_X\|_s^2 = o(1)$ as $\alpha \to 0$ and $m \to \infty$.

PROPOSITION 4.1. *Suppose that $f_X \in H_p$, $p \geq 0$. Let $\widehat{f_Y}$ be a consistent estimator of $f_Y$, that is, $\mathbb{E}\|\widehat{f_Y} - f_Y\|^2 = o(1)$ as $n \to \infty$. Consider, for $0 \leq s \leq p$, the estimator $\widehat{f_{Xs}}$ given in (2.8) with threshold $(1/m \vee \mathbb{E}\|\widehat{f_Y} - f_Y\|^2)/\alpha = o(1)$ and $\alpha = o(1)$ as $n, m \to \infty$. Then, $\mathbb{E}\|\widehat{f_{Xs}} - f_X\|_s^2 = o(1)$ as $n, m \to \infty$.*

REMARK 4.1. If we assume, in addition to the conditions of Proposition 4.1, that $m^{-1} = O(\mathbb{E}\|\widehat{f_Y} - f_Y\|^2)$ as $n \to \infty$, then we recover the result of Proposition 3.1 when $f_\epsilon$ is a priori known. In fact, in all the results below the condition $m^{-1} = O(\mathbb{E}\|\widehat{f_Y} - f_Y\|^2)$ on the sample size $m$ as $n \to \infty$, ensures that the error due to the estimation of $f_\epsilon$ is asymptotically negligible. However, in some special cases an even slower rate of $m$ is possible (see, e.g., Theorems 4.2 or 4.6).

THEOREM 4.2. *Let $f_X$ satisfy the polynomial source condition (3.2) for some $s \geq 0$ and $\beta > 0$. Consider the estimator $\widehat{f}_{X_s}$ defined in (2.8) with $\alpha = c \cdot \{(\mathbb{E}\|\widehat{f_Y} - f_Y\|^2)^{1/(\beta+1)} + m^{-1}\}$, $c > 0$. Then, for $\mathbb{E}\|\widehat{f_Y} - f_Y\|^2 \to 0$ and $m \to \infty$, we have $\mathbb{E}\|\widehat{f}_{X_s} - f_X\|_s^2 \leq C \cdot \{(\mathbb{E}\|\widehat{f_Y} - f_Y\|^2)^{\beta/(\beta+1)} + m^{-(\beta \wedge 1)}\}$, for some $C > 0$ depending only on $\rho$ given in (3.2), $\beta$ and $c$.*

REMARK 4.2. To illustrate the last result, suppose the sample size $m$ satisfies $m^{-1} = O((\mathbb{E}\|\widehat{f_Y} - f_Y\|^2)^{(\beta \vee 1)/(\beta+1)})$ as $n \to \infty$, and hence $m$ grows with a slower rate than $m^{-1} = O(\mathbb{E}\|\widehat{f_Y} - f_Y\|^2)$ (see Remark 4.1). Then, the $H_s$-risk of $\widehat{f}_{X_s}$ is bounded up to a constant by $(\mathbb{E}\|\widehat{f_Y} - f_Y\|^2)^{\beta/(\beta+1)}$, as in the case of an a priori known $f_\epsilon$ (see Theorem 3.2).

The next assertion shows that the second term given in the bound of Theorem 4.2 cannot be avoided when the samples from $f_Y$ and $f_\epsilon$ are independent. For $f \in L^2(\mathbb{R})$, let us define the class of densities

(4.2) $\quad \mathcal{D}_f^\gamma := \{g \in \mathcal{D} : \gamma |\mathcal{F}f|^2 \leq |\mathcal{F}g|^2 \leq \gamma^{-1}|\mathcal{F}f|^2\}, \qquad \gamma > 0.$

PROPOSITION 4.3. *Suppose the samples from $f_Y$ and $f_\epsilon$ are independent. Let $f \in \mathcal{D}$, and define $\mathcal{S}_f^\rho := \{g \in \mathcal{D} : \|\ell_s \cdot \mathcal{F}g \cdot (|\mathcal{F}f/\ell_s|^2)^{-\beta/2}\| \leq \rho\}$, $\rho > 0$. Then, we have $\inf_{\widehat{f_X}} \sup_{f_\epsilon \in \mathcal{D}_f^\gamma, f_X \in \mathcal{S}_f^\rho} \mathbb{E}\|\widehat{f_X} - f_X\|_s^2 \geq C \cdot m^{-(\beta \wedge 1)}$, for some $C > 0$, depending only on $f$, $\rho$ and $\gamma$.*

If $f_\epsilon$ is *ordinary smooth*, that is, (3.3) holds for some $a > 1/2$, then $f_X \in H_p$, $p > 0$ is equivalent to the polynomial source condition (3.2) with $0 \leq s < p$ and $\beta = (p-s)/(s+a)$. Thus, Theorem 4.2 implies the next assertion.

COROLLARY 4.4. *Suppose $f_\epsilon$ satisfies (3.3) for $a > 1/2$ and $f_X \in H_p$, $p > 0$. Let $\widehat{f_Y}$ given in (2.5) be constructed by using a kernel $K \in \mathcal{K}_{p+a}$ and a bandwidth $h = cn^{-1/(2(p+a)+1)}$, $c > 0$. Consider, for $0 \leq s < p$, the estimator $\widehat{f}_{X_s}$ defined in (2.8) with $\alpha = c\{n^{-2(s+a)/(2(p+a)+1)} + m^{-1}\}$, $c > 0$. Then, $\mathbb{E}\|\widehat{f}_{X_s} - f_X\|_s^2 = O(n^{-2(p-s)/(2(p+a)+1)} + m^{-(1 \wedge (p-s)/(a+s))})$ as $n, m \to \infty$.*

In case of an a priori known and ordinary smooth $f_\epsilon$, the optimal order of the $H_s$-risk over $H_p^\rho$ is $n^{-2(p-s)/(2(p+a)+1)}$ (see Remark 3.2), which together with Proposition 4.3 implies the next corollary.

COROLLARY 4.5. *Suppose the samples from $f_Y$ and $f_\epsilon$ are independent. Denote by $\mathcal{D}_a$ the set of all densities satisfying (3.3) with $a > 1/2$. Then,* $\inf_{\widehat{f_X}} \sup_{f_X \in H_p^\rho, f_\epsilon \in \mathcal{D}_a} \mathbb{E}\|\widehat{f_X} - f_X\|_s^2 \geq C\{n^{-2(p-s)/(2(p+a)+1)} + m^{-(1 \wedge (p-s)/(a+s))}\}$.

REMARK 4.3. If the samples from $f_Y$ and $f_\epsilon$ are independent, then due to Corollaries 4.4 and 4.5 the order of the smallest $m$ for archiving the same convergence rate as in the case of an a priori known $f_\epsilon$ (Corollary 3.3) is given by $m^{-1} = O(n^{-2[(p-s)\vee(a+s)]/[2(p+a)+1]})$. We shall emphasize the interesting ambiguous influences of the parameters $p$ and $a$ characterizing the smoothness of $f_X$ and $f_\epsilon$, respectively. If in case of $(p-s) < (a+s)$ the value of $a$ decreases or the value of $p$ increases, then the estimation of $f_\epsilon$ is still negligible given a relative to $n$ slower necessary rate of $m$. While in the case of $(p-s) > (a+s)$ a decreasing value of $a$ or an increasing value of $p$ leads to a relative to $n$ faster necessary rate of $m$. However, in both cases a decreasing value of $a$ or an increasing value of $p$ implies a faster optimal rate of convergence of the estimator $\widehat{f_X}_s$.

THEOREM 4.6. *Let $f_X$ satisfy the logarithmic source condition (3.4), for some $s \geq 0$ and $\beta > 0$. Consider the estimator $\widehat{f_X}_s$ defined in (2.8) by using a threshold $\alpha = c\{(\mathbb{E}\|\widehat{f_Y} - f_Y\|^2)^{1/2} + m^{-1/2}\}$, $c > 0$. Then, for $\mathbb{E}\|\widehat{f_Y} - f_Y\|^2 \to 0$ and $m \to \infty$, we have $\mathbb{E}\|\widehat{f_X}_s - f_X\|_s^2 \leq C\{|\log(\mathbb{E}\|\widehat{f_Y} - f_Y\|^2)|^{-\beta} + (\log m)^{-\beta}\}$, for some $C > 0$ depending only on $\rho$ given in (3.4), $\beta$ and $c$.*

REMARK 4.4. Assume that, for some $\nu > 0$, the sample size $m$ satisfies $m^{-1} = O((\mathbb{E}\|\widehat{f_Y} - f_Y\|^2)^\nu)$ as $n \to \infty$, and hence $m$ may grow with a fare slower rate than implied by the condition $m^{-1} = O(\mathbb{E}\|\widehat{f_Y} - f_Y\|^2)$ (compare Remark 4.1). Then, as in the case of an a priori known $f_\epsilon$ (see Theorem 3.4), the $H_s$-risk of $\widehat{f_X}_s$ is bounded by $C|\log(\mathbb{E}\|\widehat{f_Y} - f_Y\|^2)|^{-\beta}$, for some $C > 0$. Note that the influence of the parameter $\nu$ is hidden in the constant $C$.

The next assertion states that the second term given in the bound of Theorem 4.6 cannot be avoided.

PROPOSITION 4.7. *Suppose the samples from $f_Y$ and $f_\epsilon$ are independent. Let $f \in \mathcal{D}$, and define $\mathcal{S}_f^\rho := \{g \in \mathcal{D} : \|\ell_s \cdot \mathcal{F}g \cdot |\log(|\mathcal{F}f/\ell_s|^2)|^{\beta/2}\| \leq \rho\}$, $\rho > 0$. Then, we have $\inf_{\widehat{f_X}} \sup_{f_\epsilon \in \mathcal{D}_f^\gamma, f_X \in \mathcal{S}_f^\rho} \mathbb{E}\|\widehat{f_X} - f_X\|_s^2 \geq C(\log m)^{-\beta}$, for some $C > 0$, depending only on $f$, $\rho$ and $\gamma$.*

Assume that $f_\epsilon$ is supersmooth, that is, (3.5) holds for $a > 0$. Then, $f_X \in H_p$, $p > 0$, is equivalent to the logarithmic source condition (3.4) with $0 \leq s < p$ and $\beta = (p-s)/a$. Thus, Theorem 4.6 implies the next assertion.

COROLLARY 4.8. *Suppose $f_\epsilon$ satisfies (3.5), for $a > 0$ and $f_X \in H_p$, $p > 0$. Let $\widehat{f_Y}$ defined in (2.5) be constructed by using a kernel $K \in \mathcal{K}_r$ [see (2.6)] and a bandwidth $h = cn^{-1/(2r+1)}$, $c, r > 0$. Consider, for $0 \leq s < p$, the estimator $\widehat{f_X}_s$ defined in (2.8) with $\alpha = c\{n^{-r/(2r+1)} + m^{-1/2}\}$, $c > 0$. Then, $\mathbb{E}\|\widehat{f_X}_s - f_X\|_s^2 = O((\log n)^{-(p-s)/a} + (\log m)^{-(p-s)/a})$ as $n, m \to \infty$.*

In case of an a priori known and supersmooth $f_\epsilon$, the optimal order of the $H_s$-risk over $H_p^\rho$ is $(\log n)^{-(p-s)/a}$ (see Remark 3.3), which together with Proposition 4.7 leads to the next assertion.

COROLLARY 4.9. *Suppose the samples from $f_Y$ and $f_\epsilon$ are independent. Denote by $\mathcal{D}_a$ the set of all densities satisfying (3.5) with $a > 0$. Then, $\inf_{\widehat{f_X}} \sup_{f_X \in H_p^\rho, f_\epsilon \in \mathcal{D}_a} \mathbb{E}\|\widehat{f_X} - f_X\|_s^2 \geq C\{(\log n)^{-(p-s)/a} + (\log m)^{-(p-s)/a}\}.$*

REMARK 4.5. If we assume $m^{-1} = O(n^{-\nu})$, for some $\nu > 0$, then the order in the last result simplifies to $(\log n)^{-(p-s)/a}$ and hence, equals the optimal order for known $f_\epsilon$ (see Corollary 3.5). Therefore, if the samples from $f_Y$ and $f_\epsilon$ are independent, then from Corollary 4.8 and 4.9 it follows that the error due to the estimation of $f_\epsilon$ is asymptotically negligible if and only if the sample size $m$ grows as some power of $n$. In contrast to the situation in Corollary 4.4 and 4.5, if $f_\epsilon$ is supersmooth, that is, (3.5) holds for $a > 0$, and $f_X \in H_p$, $p > 0$, then the influence of the parameters $p$ and $a$ is not ambiguous. A decreasing value of $a$ or an increasing value of $p$ implies a faster optimal rate of convergence of the estimator $\widehat{f_X}_s$, and the relative to $n$ necessary rate of $m$ is not affected. Note that the estimator is adaptive as in a case of known supersmooth error density (see Remark 3.3). We shall stress that the estimation of $f_\epsilon$ has no influence on the order of the $H_s$-risk of $\widehat{f_X}_s$, as long as the sample size $m$ grows as fast as some power of $n$. However, the influence is clearly hidden in the constant of the order symbol.

THEOREM 4.10. *Let $f_X$ satisfy the general source condition (3.6) for some concave index function $\kappa$ and $s \geq 0$. Denote by $\Phi$ and $\omega$ the inverse function of $\kappa$ and $\omega^{-1}(t) := t\Phi(t)$, respectively. Consider $\widehat{f_X}_s$ defined in (2.8) with $\alpha = c\{\mathbb{E}\|\widehat{f_Y} - f_Y\|^2/\omega(\mathbb{E}\|\widehat{f_Y} - f_Y\|^2) + 1/m\}$, $c > 0$. Then, we have $\mathbb{E}\|\widehat{f_X}_s - f_X\|_s^2 \leq C\{\omega(\mathbb{E}\|\widehat{f_Y} - f_Y\|^2) + \kappa(1/m)\}$, as $\mathbb{E}\|\widehat{f_Y} - f_Y\|^2 \to 0$ and $m \to \infty$, for some $C > 0$, depending only on $\rho$ given in (3.6) and c.*

REMARK 4.6. Assume that $m^{-1} = O(\mathbb{E}\|\widehat{f_Y} - f_Y\|^2)$ as $n \to \infty$, then the $H_s$-risk of $\widehat{f}_{X_s}$ is bounded up to a constant by $\omega(\mathbb{E}\|\widehat{f_Y} - f_Y\|^2)$ as in case of an a priori known $f_\epsilon$ (see Theorem 3.6). Thus, the general source condition supposing $m^{-1} = O(\mathbb{E}\|\widehat{f_Y} - f_Y\|^2)$ is sufficient to ensure that the estimation of the noise density is asymptotically negligible.

PROPOSITION 4.11. *Let the samples from $f_Y$ and $f_\epsilon$ be independent and $f \in \mathcal{D}$. Define $\mathcal{S}_f^\rho := \{g \in \mathcal{D} : \|\ell_s \cdot \mathcal{F}g \cdot \kappa(|\mathcal{F}f/\ell_s|^2)|^{-1/2}\| \leq \rho\}$, $\rho > 0$. Then, we have $\inf_{\widehat{f_X}} \sup_{f_\epsilon \in \mathcal{D}_f^\gamma, f_X \in \mathcal{S}_f^\rho} \mathbb{E}\|\widehat{f_X} - f_X\|_s^2 \geq C \cdot \kappa(1/m)$, for some $C > 0$ depending only on $f$, $\rho$ and $\gamma$.*

REMARK 4.7. Due to Proposition 4.11 in the case of independent samples from $f_Y$ and $f_\epsilon$, the term $\kappa(1/m)$ given in the bound of Theorem 4.10 cannot be avoided. It follows that our estimator $\widehat{f}_{X_s}$ attains the minimax optimal order over $\mathcal{S}_{f_\epsilon}^\rho$ when $\omega(\mathbb{E}\|\widehat{f_Y} - f_Y\|^2)$ is the optimal order for known $f_\epsilon$ (compare Remark 3.4).

## APPENDIX

PROOF OF PROPOSITION 3.1. The proof is based on the decomposition (3.1), where $\alpha^{-1} \geq \sup_{t \in \mathbb{R}^+} t^{-1} \mathbb{1}\{t \geq \alpha\}$ is used to obtain the first term on the right-hand side. If $f_X \in H_p$, $p \geq s \geq 0$, then by making use of the relation $\|f_{X_s}^\alpha - f_X\|_s^2 = \|\mathbb{1}\{|\mathcal{F}f_\epsilon/\ell_s|^2 < \alpha\} \cdot \ell_s \cdot \mathcal{F}f_X\|^2 \leq \|\ell_s \cdot \mathcal{F}f_X\|^2 \leq \|f_X\|_p^2 < \infty$, the second term satisfies $\|f_{X_s}^\alpha - f_X\|_s^2 = o(1)$, as $\alpha \to 0$, due to Lebesgue's dominated convergence theorem. Therefore, the conditions on $\alpha$ ensure the convergence to zero of the two terms on the right-hand side in (3.1) as $n$ increases, which gives the result. □

Assuming $f_\epsilon$ is known, the next lemma summarizes the essential bounds of the regularization bias depending on the polynomial, logarithmic or general source condition.

LEMMA A.1. *Let $w : \mathbb{R} \to [1, \infty)$ be an arbitrary weight function. Suppose there exists $\beta > 0$ such that:*

(i) $\rho := \|w \cdot \mathcal{F}f_X \cdot (|\mathcal{F}f_\epsilon|^2/w^2)^{-\beta/2}\| < \infty$ *is satisfied, then*

(A.1) $\qquad \|w \cdot \mathcal{F}f_X \cdot \mathbb{1}\{|\mathcal{F}f_\epsilon|^2/w^2 < \alpha\}\|^2 \leq \alpha^\beta \cdot \rho^2;$

(ii) $\rho := \|w \cdot \mathcal{F}f_X \cdot |\log(|\mathcal{F}f_\epsilon|^2/w^2)|^{\beta/2}\| < \infty$ *is satisfied, then*

(A.2) $\qquad \|w \cdot \mathcal{F}f_X \cdot \mathbb{1}\{|\mathcal{F}f_\epsilon|^2/w^2 < \alpha\}\|^2 \leq C_\beta \cdot (-\log \alpha)^{-\beta} \cdot \rho^2;$

(iii) $\rho := \|w \cdot \mathcal{F} f_X \cdot |\kappa(|\mathcal{F} f_\epsilon|^2/w^2)|^{-1/2}\| < \infty$ is satisfied and assume that the index function $\kappa$ is concave, then

$$(A.3) \qquad \|w \cdot \mathcal{F} f_X \cdot \mathbb{1}\{|\mathcal{F} f_\epsilon|^2/w^2 < \alpha\}\|^2 \leq C_\kappa \cdot \kappa(\alpha) \cdot \rho^2;$$

where $C_\beta$, $C_\kappa$ are positive constants depending only on $\beta$ and $\kappa$, respectively.

PROOF. Denote $\psi_\alpha := \mathcal{F} f_X \mathbb{1}\{|\mathcal{F} f_\epsilon/w|^2 < \alpha\}$. Under the assumption (i) we have $\|w \cdot \psi_\alpha\|^2 \leq \sup_{t \in \mathbb{R}^+} t^\beta \mathbb{1}\{t < \alpha\} \cdot \rho^2$, which implies (A.1).

The proof of (A.2) is partially motivated by techniques used in Nair, Pereverzev and Tautenhahn [29]. Let $\kappa_\beta(t) := |\log(t)|^{-\beta}$, $t \in (0,1)$ and $\phi_\beta(t) := \kappa_\beta^{1/2}(|[\mathcal{F} f_\epsilon](t)/w(t)|^2)$, $t \in \mathbb{R}$, then for all $t \in \mathbb{R}$ we have $\phi_\beta(0) \geq \phi_\beta(t) > 0$. Under assumption (ii), which may be rewritten as $\rho = \|w \cdot \mathcal{F} f_X / \phi_\beta\| < \infty$, we obtain

$$(A.4) \qquad \|w \cdot \psi_\alpha\|^2 = \int_\mathbb{R} w(t) \psi_\alpha(t) \phi_\beta(t) \frac{w(t) \overline{[\mathcal{F} f_X](t)}}{\phi_\beta(t)} \, dt \leq \|w \cdot \psi_\alpha \cdot \phi_\beta\| \cdot \rho$$

due to the Cauchy–Schwarz inequality. From (A.4) we conclude

$$(A.5) \qquad \|\mathcal{F} f_\epsilon \cdot \psi_\alpha\|^2 = \|\mathcal{F} f_\epsilon \cdot \mathbb{1}\{|\mathcal{F} f_\epsilon/w|^2 < \alpha\} \cdot \psi_\alpha\|^2 \leq \alpha \cdot \|w \cdot \psi_\alpha \cdot \phi_\beta\| \cdot \rho,$$

since $\alpha \geq \sup_{t \in \mathbb{R}^+} t \cdot \mathbb{1}\{t < \alpha\}$. Let $\Phi_\beta$ be the inverse function of $\kappa_\beta$, then $\Phi_\beta(s) = e^{-s^{-1/\beta}}$, $s > 0$, which is convex on the interval $(0, c_\beta^2]$ with $c_\beta^2 = (1+\beta)^{-\beta}$. Define $\gamma_\beta^2 = c_\beta^2/\phi_\beta^2(0) \wedge 1$. Therefore, Jensen's inequality implies

$$\Phi_\beta\left(\frac{\gamma_\beta^2 \cdot \|w \cdot \psi_\alpha \cdot \phi_\beta\|^2}{\|w \cdot \psi_\alpha\|^2}\right) \leq \frac{\int_\mathbb{R} \Phi_\beta(\gamma_\beta^2 \cdot \phi_\beta^2(t)) \cdot w^2(t) \cdot \psi_\alpha^2(t) \, dt}{\int_\mathbb{R} w^2(t) \cdot \psi_\alpha^2(t) \, dt},$$

which together with $\Phi_\beta(\gamma_\beta^2 \cdot \phi_\beta^2(t)) \leq \Phi_\beta(\phi_\beta^2(t)) = |[\mathcal{F} f_\epsilon](t)|^2/w^2(t)$ gives

$$(A.6) \qquad \Phi_\beta\left(\frac{\gamma_\beta^2 \cdot \|w \cdot \psi_\alpha \cdot \phi_\beta\|^2}{\|w \cdot \psi_\alpha\|^2}\right) \leq \frac{\int_\mathbb{R} |[\mathcal{F} f_\epsilon](t)|^2 \cdot \psi_\alpha^2(t) \, dt}{\|w \cdot \psi_\alpha\|^2} = \frac{\|\mathcal{F} f_\epsilon \cdot \psi_\alpha\|^2}{\|w \cdot \psi_\alpha\|^2}.$$

In order to combine the three estimates (A.4), (A.5) and (A.6), let us introduce a new function $\Psi_\beta$ by $\Psi_\beta(t) := \Phi_\beta(t^2)/t^2$. Since $\Phi_\beta$ is convex, we conclude that $\Psi_\beta$ is monotonically increasing on the interval $(0, c_\beta]$. Hence, by (A.4), which may be rewritten as $\|w \cdot \psi_\alpha \cdot \phi_\beta\|^{1/2}/\rho^{1/2} \leq \|w \cdot \psi_\alpha \cdot \phi_\beta\|/\|w \cdot \psi_\alpha\|$ ($\leq \phi_\beta(0)$), the monotonicity of $\Psi_\beta$ and (A.6),

$$\Psi_\beta\left(\frac{\gamma_\beta \cdot \|w \cdot \psi_\alpha \cdot \phi_\beta\|^{1/2}}{\rho^{1/2}}\right) \leq \Psi_\beta\left(\frac{\gamma_\beta \cdot \|w \cdot \psi_\alpha \cdot \phi_\beta\|}{\|w \cdot \psi_\alpha\|}\right) \leq \frac{\|\mathcal{F} f_\epsilon \cdot \psi_\alpha\|^2}{\gamma_\beta^2 \cdot \|w \cdot \psi_\alpha \cdot \phi_\beta\|^2}.$$

Multiplying by $\gamma_\beta^2 \cdot \|w \cdot \psi_\alpha \cdot \phi_\beta\|/\rho$ and exploiting (A.5) yields

$$(A.7) \qquad \Phi_\beta\left(\frac{\gamma_\beta^2 \cdot \|w \cdot \psi_\alpha \cdot \phi_\beta\|}{\rho}\right) \leq \frac{\|\overline{\mathcal{F} f_\epsilon} \cdot \psi_\alpha\|^2}{\rho \cdot \|w \cdot \psi_\alpha \cdot \phi_\beta\|} \leq \alpha.$$

Since $\Phi_\beta^{-1}(s) = |\ln(s)|^{-\beta}$, we obtain (A.2) by combining (A.4) and (A.7).

The proof of (A.3) follows line by line the proof of (A.2) using the concave index function $\kappa$ and its convex inverse function $\Phi$, rather than $\kappa_\beta$ and $\Phi_\beta$. □

PROOF OF THEOREM 3.2. The proof is based on the decomposition (3.1). The polynomial source condition (3.2) equals assumption (i) in Lemma A.1 with $w \equiv \ell_s$, therefore from (A.1) we obtain $\|f_{X\,s}^\alpha - f_X\|_s^2 \leq \alpha^\beta \cdot \rho^2$. Balancing the two terms on the right-hand side in (3.1) then gives the result. □

PROOF OF COROLLARY 3.3. Under the conditions of the corollary, we have $f_Y \in H_{p+a}$ and, hence $\mathbb{E}\|\widehat{f_Y} - f_Y\|^2 = O(n^{-2(p+a)/(2(p+a)+1)})$. Moreover, the polynomial source condition (3.2) is satisfied with $\beta = (p-s)/(a+s)$. Therefore, the result follows from Theorem 3.2. □

PROOF OF THEOREM 3.4. The proof is similar to the proof of Theorem 3.2, but uses (A.2) in Lemma A.1 with $w \equiv \ell_s$ rather than (A.1). The conditions of the theorem then provide $\mathbb{E}\|\widetilde{f_X}_s - f_X\|_s^2 \leq C(\mathbb{E}\|\widehat{f_Y} - f_Y\|^2)^{1/2} + C|\log(\mathbb{E}\|\widehat{f_Y} - f_Y\|^2)|^{-\beta}$, for some constant $C > 0$, depending only on $\rho$ given in (3.4), $\beta$ and $c$, which implies the result. □

PROOF OF COROLLARY 3.5. Under the conditions of the corollary, we have $f_Y \in H_r$ and, hence $\mathbb{E}\|\widehat{f_Y} - f_Y\|^2 = O(n^{-2r/(2r+1)})$. Moreover, the logarithmic source condition (3.4) is satisfied with $\beta = (p-s)/a$. Therefore, the result follows from Theorem 3.4. □

PROOF OF THEOREM 3.6. The proof is similar to the proof of Theorem 3.2, but uses (A.3) in Lemma A.1 with $w \equiv \ell_s$ rather than (A.1). The condition on $\alpha$ which may be rewritten as $c \cdot \mathbb{E}\|\widehat{f_Y} - f_Y\|^2 = \alpha \cdot \kappa(\alpha)$ then ensures the balance of the two terms in (3.1). The result follows by making use of the relation $\omega(c \cdot \delta) \leq (c \vee 1) \cdot \omega(\delta)$ (Mair and Ruymgaart [26], Remark 3.7). □

LEMMA A.2. *Suppose $w:\mathbb{R} \to [1,\infty)$ is an arbitrary weight function, $\kappa$ is a concave index function and $\widehat{\mathcal{F}f_\epsilon}$ is the estimator defined in (2.7). Then, for all $\gamma \geq 0$ and $t \in \mathbb{R}$, we have*

$$\mathbb{E}|[\widehat{\mathcal{F}f_\epsilon}](t)/w(t) - [\mathcal{F}f_\epsilon](t)/w(t)|^{2\gamma} \tag{A.8}$$
$$\leq C(\gamma) \cdot m^{-\gamma},$$

$$\mathbb{E}\left[\mathbb{1}\{|[\widehat{\mathcal{F}f_\epsilon}](t)/w(t)|^2 \geq \alpha\} \cdot \frac{|[\widehat{\mathcal{F}f_\epsilon}](t) - [\mathcal{F}f_\epsilon](t)|^2}{|[\widehat{\mathcal{F}f_\epsilon}](t)|^2}\right]$$

(A.9)
$$\leq \frac{C(\gamma)}{|[\mathcal{F}f_\epsilon](t)/w(t)|^{2\gamma}} \cdot \left\{ \frac{1}{\alpha \cdot m^{1+\gamma}} + \frac{1}{(m \cdot \alpha)^{1-\gamma\wedge 1} \cdot m^{\gamma\wedge 1}} \right\},$$

(A.10)
$$\mathbb{E}\left[\mathbb{1}\{|[\widehat{\mathcal{F}f_\epsilon}](t)/w(t)|^2 \geq \alpha\} \cdot \frac{|[\widehat{\mathcal{F}f_\epsilon}](t) - [\mathcal{F}f_\epsilon](t)|^2}{|[\widehat{\mathcal{F}f_\epsilon}](t)|^2}\right]$$
$$\leq \frac{C(\gamma)}{\kappa(|[\mathcal{F}f_\epsilon](t)/w(t)|^2)} \cdot \left\{ \frac{\kappa(1/m)}{\alpha \cdot m} + \kappa(1/m) \right\},$$

where $C$ and $C(\gamma)$ depending only on $\gamma$ are positive constants.

PROOF. Let $\gamma \geq 0$ and $t \in \mathbb{R}$. Define $Z_j := \{(2\pi)^{-1/2} e^{-it\epsilon_j} - [\mathcal{F}f_\epsilon](t)\}/w(t)$, $j = 1, \ldots, m$, then $Z_1, \ldots, Z_m$ are i.i.d. random variables with mean zero, and $|Z_j|^{2\gamma} \leq K$ for some positive constant $K$. Therefore, applying Theorem 2.10 in Petrov [35], we obtain (A.8) for $\gamma \geq 1$, while for $\gamma \in (0, 1)$ the estimate follows from Lyapunov's inequality.

Proof of (A.9). Consider, for $\gamma \geq 0$ and $t \in \mathbb{R}$, the elementary inequality

(A.11) $$1 \leq 2^{2\gamma} \cdot \left\{ \frac{|[\widehat{\mathcal{F}f_\epsilon}](t)/w(t) - [\mathcal{F}f_\epsilon](t)/w(t)|^{2\gamma}}{|[\mathcal{F}f_\epsilon](t)/w(t)|^{2\gamma}} + \frac{|[\widehat{\mathcal{F}f_\epsilon}](t)/w(t)|^{2\gamma}}{|[\mathcal{F}f_\epsilon](t)/w(t)|^{2\gamma}} \right\},$$

which together with $|[\widehat{\mathcal{F}f_\epsilon}](t)/w(t)| \leq 1$, for all $t \in \mathbb{R}$, implies

$$\mathbb{E}\left[\mathbb{1}\{|[\widehat{\mathcal{F}f_\epsilon}](t)/w(t)|^2 \geq \alpha\} \cdot \frac{|[\widehat{\mathcal{F}f_\epsilon}](t) - [\mathcal{F}f_\epsilon](t)|^2}{|[\widehat{\mathcal{F}f_\epsilon}](t)|^2}\right]$$
$$\leq \frac{2^{2\gamma}}{|[\mathcal{F}f_\epsilon](t)/w(t)|^{2\gamma}} \cdot \left\{ \frac{\mathbb{E}|[\widehat{\mathcal{F}f_\epsilon}](t)/w(t) - [\mathcal{F}f_\epsilon](t)/w(t)|^{2(1+\gamma)}}{\alpha} \right.$$
$$\left. + \frac{\mathbb{E}|[\widehat{\mathcal{F}f_\epsilon}](t)/w(t) - [\mathcal{F}f_\epsilon](t)/w(t)|^2}{\alpha^{1-\gamma\wedge 1}} \right\}$$

and by using (A.8) we obtain the estimate (A.9).

Proof of (A.9). If $|[\mathcal{F}f_\epsilon](t)/w(t)|^2 \leq 1/m$, then we obtain (A.10) by using (A.8) with $\gamma = 1$ together with $\kappa(|[\mathcal{F}f_\epsilon](t)/w(t)|^2) \leq \kappa(1/m)$. Since $\kappa$ is concave, we conclude that $g(t) = \kappa(t^2)/t^2$ is monotonically decreasing. Hence, if $|[\mathcal{F}f_\epsilon](t)/w(t)|^2 \geq 1/m$, then due to the monotonicity of $g$ we have $\kappa(|[\mathcal{F}f_\epsilon](t)/w(t)|^2)|[\mathcal{F}f_\epsilon](t)/w(t)|^{-2} \leq m\kappa(m^{-1})$, which together with inequality (A.11), for $\gamma = 1$, yields

$$\mathbb{E}\left[\mathbb{1}\{|[\widehat{\mathcal{F}f_\epsilon}](t)/w(t)|^2 \geq \alpha\} \cdot \frac{|[\widehat{\mathcal{F}f_\epsilon}](t) - [\mathcal{F}f_\epsilon](t)|^2}{|[\widehat{\mathcal{F}f_\epsilon}](t)|^2}\right]$$
$$\leq \frac{2^4 m\kappa(m^{-1})}{\kappa(|[\mathcal{F}f_\epsilon](t)/w(t)|^2)} \cdot \left\{ \frac{\mathbb{E}|[\widehat{\mathcal{F}f_\epsilon}](t)/w(t) - [\mathcal{F}f_\epsilon](t)/w(t)|^4}{\alpha} \right.$$

$$+ \mathbb{E}|[\widehat{\mathcal{F}f_\epsilon}](t)/w(t) - [\mathcal{F}f_\epsilon](t)/w(t)|^2 \Big\}$$

and by using (A.8) we obtain the estimate (A.10). □

PROOF OF PROPOSITION 4.1. The proof is based on the decomposition (4.1). Due to (A.9) in Lemma A.2, we show below the bound

$$\text{(A.12)} \quad \begin{aligned} \mathbb{E}\|\widehat{f_X}_s - \widehat{f_X^\alpha}_s\|_s^2 &\leq \pi^{-1}\alpha^{-1} \cdot \mathbb{E}\|\widehat{f_Y} - f_Y\|^2 \\ &\quad + 2C(0) \cdot \|f_X\|_s^2 \cdot \alpha^{-1} \cdot m^{-1}, \end{aligned}$$

while from Lebesgue's dominated convergence theorem and (A.8) in Lemma A.2, we conclude

$$\text{(A.13)} \quad \mathbb{E}\|\widehat{f_X^\alpha}_s - f_X\|_s^2 = o(1) \qquad \text{as } \alpha \to 0 \text{ and } m \to \infty.$$

Therefore, the conditions on $\alpha$ ensure the convergence to zero of the two terms on the right-hand side in (4.1) as $n$ and $m$ tend to $\infty$, which gives the result.

Proof of (A.12). Using $\alpha^{-1} \geq \sup_{t \in \mathbb{R}^+} t^{-1}\mathbb{1}\{t \geq \alpha\}$, we have

$$\text{(A.14)} \quad \begin{aligned} &\mathbb{E}\|\widehat{f_X}_s - \widehat{f_X^\alpha}_s\|_s^2 \\ &\leq \pi^{-1}\alpha^{-1} \cdot \mathbb{E}\|\mathcal{F}\widehat{f_Y} - \mathcal{F}f_Y\|^2 \\ &\quad + 2\bigg\|\bigg\{\mathbb{E}\bigg[\mathbb{1}\{|\widehat{\mathcal{F}f_\epsilon}/\ell^s|^2 \geq \alpha\} \cdot \frac{|\widehat{\mathcal{F}f_\epsilon}/\ell_s - \mathcal{F}f_\epsilon/\ell_s|^2}{|\widehat{\mathcal{F}f_\epsilon}/\ell_s|^2}\bigg]\bigg\}^{1/2} \\ &\qquad\qquad \times \ell_s \cdot \mathcal{F}f_X\bigg\|^2 \end{aligned}$$

and hence $\|\ell_s \cdot \mathcal{F}f_X\| = \|f_X\|_s \leq \|f_X\|_p < \infty$, together with (A.9) in Lemma A.2 with $w = \ell_s$ and $\gamma = 0$, implies (A.12).

Proof of (A.13). If $f_X \in H_p$, $p \geq s \geq 0$, then by making use of the relation $\mathbb{E}\|\widehat{f_X^\alpha}_s - f_X\|_s^2 = \|\mathbb{E}\mathbb{1}\{|\widehat{\mathcal{F}f_\epsilon}/\ell_s|^2 < \alpha\} \cdot \ell_s \cdot \mathcal{F}f_X\|^2 \leq \|\ell_s \cdot \mathcal{F}f_X\|^2 \leq \|f_X\|_p^2 < \infty$ the result follows due to Lebesgue's dominated convergence theorem from $\mathbb{E}\mathbb{1}\{|[\widehat{\mathcal{F}f_\epsilon}](t)/\ell_s(t)|^2 < \alpha\} \to 0$ as $\alpha \to 0$ and $m \to \infty$, that can be realized as follows. For all $\alpha \leq \alpha_0$, we have $|[\mathcal{F}f_\epsilon](t)| \geq 2\alpha^{1/2}\ell_s(t)$ and, hence $\mathbb{E}\mathbb{1}\{|[\widehat{\mathcal{F}f_\epsilon}](t)/\ell_s(t)|^2 < \alpha\} \leq P(|[\widehat{\mathcal{F}f_\epsilon}](t) - [\mathcal{F}f_\epsilon](t)| > |[\mathcal{F}f_\epsilon](t)|/2)$. Therefore, from Chebyshev's inequality and (A.8) in Lemma A.2 with $w \equiv 1$ and $\gamma = 1$, we obtain (A.13). □

The next lemma summarizes the essential bounds of the "regularization error" depending on the polynomial, logarithmic or general source condition.

LEMMA A.3. *Let $w:\mathbb{R} \to [1,\infty)$ be an arbitrary weight function, and let $\widehat{\mathcal{F}f_\epsilon}$ be the estimator defined in (2.7). Suppose there exists $\beta > 0$ such that:*

(i) $\rho := \|w \cdot \mathcal{F}f_X \cdot (|\mathcal{F}f_\epsilon|^2/w^2)^{-\beta/2}\| < \infty$ *is satisfied, then*

$$(A.15) \quad \mathbb{E}\|w \cdot \mathcal{F}f_X \cdot \mathbb{1}\{|\widehat{\mathcal{F}f_\epsilon}/w|^2 < \alpha\}\|^2 \leq C_\beta\{\alpha^\beta + m^{-\beta}\}\rho^2;$$

(ii) $\rho := \|w \cdot \mathcal{F}f_X \cdot |\log(|\mathcal{F}f_\epsilon/w|^2)|^{\beta/2}\| < \infty$ *is satisfied, then*

$$(A.16) \quad \mathbb{E}\|w \cdot \mathcal{F}f_X \cdot \mathbb{1}\{|\widehat{\mathcal{F}f_\epsilon}/w|^2 < \alpha\}\|^2 \leq C_\beta |\log(C_\beta\{\alpha + m^{-1}\})|^{-\beta}\rho^2;$$

(iii) $\rho := \|w \cdot \mathcal{F}f_X \cdot |\kappa(|\mathcal{F}f_\epsilon/w|^2)|^{-1/2}\| < \infty$, *and assume that the index function $\kappa$ is concave, then*

$$(A.17) \quad \mathbb{E}\|w \cdot \mathcal{F}f_X \cdot \mathbb{1}\{|\widehat{\mathcal{F}f_\epsilon}/w|^2 < \alpha\}\|^2 \leq C_\kappa \cdot \kappa(C_\kappa\{\alpha + m^{-1}\}) \cdot \rho^2;$$

*where $C_\beta$, $C_\kappa$ are positive constants depending only on $\beta$ and $\kappa$, respectively.*

PROOF. Denote $\widehat{\psi}_\alpha := \mathcal{F}f_X \cdot \mathbb{1}\{|\widehat{\mathcal{F}f_\epsilon}/w|^2 < \alpha\}$. Then, using the inequality (A.11) together with $\alpha^\gamma \geq \sup_{t \in \mathbb{R}^+} t^\gamma \mathbb{1}\{t < \alpha\}$, for all $\gamma > 0$, we have

$$\|w \cdot \widehat{\psi}_\alpha\|^2 \leq 2^{2\beta}\{\alpha^\beta \cdot \rho^2 + \|w \cdot \mathcal{F}f_X \cdot |\mathcal{F}f_\epsilon/w|^{-\beta} \cdot |\widehat{\mathcal{F}f_\epsilon}/w - \mathcal{F}f_\epsilon/w|^\beta\|^2\}.$$

Therefore, using (A.8) in Lemma A.2, we obtain the bound (A.15).

The proof of (A.16) follows along the same lines as the proof of (A.2) in Lemma A.1. Consider the functions $\kappa_\beta$, $\phi_\beta$ and $\Phi_\beta$ defined in the proof of (A.2) in Lemma A.1, then in analogy to (A.4), we bound

$$(A.18) \quad \|w \cdot \widehat{\psi}_\alpha\|^2 \leq \|w \cdot \widehat{\psi}_\alpha \cdot \phi_\beta\| \cdot \rho,$$

which implies

$$(A.19) \quad \mathbb{E}\|\widehat{\psi}_\alpha\|^2 \leq (\mathbb{E}\|w \cdot \widehat{\psi}_\alpha \cdot \phi_\beta\|^2)^{1/2} \cdot \rho.$$

Moreover, following the steps in (A.5) together with (A.18), we have

$$(A.20) \quad \|\widehat{\mathcal{F}f_\epsilon} \cdot \widehat{\psi}_\alpha\|^2 \leq \alpha \cdot \|w \cdot \widehat{\psi}_\alpha \cdot \phi_\beta\| \cdot \rho.$$

Therefore, applying the triangular inequality together with (A.20), we obtain

$$\mathbb{E}\|\mathcal{F}f_\epsilon \cdot \widehat{\psi}_\alpha\|^2 \leq 2\mathbb{E}\|w \cdot |\mathcal{F}f_\epsilon/w - \widehat{\mathcal{F}f_\epsilon}/w| \cdot \widehat{\psi}_\alpha\|^2 + 2\alpha(\mathbb{E}\|w \cdot \widehat{\psi}_\alpha \cdot \phi_\beta\|^2)^{1/2}\rho.$$

By applying the Cauchy–Schwarz inequality and then (A.8) in Lemma A.2, we bound the first term by $C(\beta) \cdot m^{-1} \cdot \int (\mathbb{E}\mathbb{1}\{|[\widehat{\mathcal{F}f_\epsilon}](t)/w(t)|^2 < \alpha\})^{1/2} \cdot w^2(t) \cdot |[\mathcal{F}f_X](t)|^2\,dt$, and using once again the Cauchy–Schwarz inequality,

$$(A.21) \quad \mathbb{E}\|\mathcal{F}f_\epsilon \cdot \widehat{\psi}_\alpha\|^2 \leq 2\left\{\frac{C(\beta)}{m} + \alpha\right\} \cdot (\mathbb{E}\|w \cdot \widehat{\psi}_\alpha \cdot \phi_\beta\|^2)^{1/2} \cdot \rho.$$

In analogy to (A.6), by applying the convex function $\Phi_\beta$, we obtain

$$\text{(A.22)} \qquad \Phi_\beta\left(\frac{\gamma_\beta^2 \cdot \mathbb{E}\|w \cdot \widehat{\psi}_\alpha \cdot \phi_\beta\|^2}{\mathbb{E}\|w \cdot \widehat{\psi}_\alpha\|^2}\right) \leq \frac{\mathbb{E}\|\mathcal{F}f_\epsilon \cdot \widehat{\psi}_\alpha\|^2}{\mathbb{E}\|w \cdot \widehat{\psi}_\alpha\|^2}.$$

Combining the three bounds (A.19), (A.21) and (A.22), as in (A.7), implies

$$\text{(A.23)} \qquad \Phi_\beta\left(\frac{\gamma_\beta^2 \cdot (\mathbb{E}\|w \cdot \widehat{\psi}_\alpha \cdot \phi_\beta\|^2)^{1/2}}{\rho}\right) \leq 2\left\{\frac{C(\beta)}{m} + \alpha\right\}.$$

We obtain the second bound (A.16) by combining (A.19) and (A.23).

The proof of (A.17) follows line by line the proof of (A.16) using the functions $\kappa$ and $\Phi$ rather than $\kappa_\beta$ and $\Phi_\beta$. $\square$

The next lemma generalizes Theorem 3.1 given in Neumann [31] by providing a lower bound for the MISE under a general source condition, which requires for $f \in L^2(\mathbb{R})$ and index function $\kappa$ the following definitions:

$$\text{(A.24)} \quad \begin{aligned} \mathcal{M}_f^\rho &:= \{g \in \mathcal{D} : \|\mathcal{F}g \cdot |\kappa(|\mathcal{F}f|^2)|^{-1/2}\| \leq \rho\}, \qquad \rho > 0, \\ \Delta_f^m(t) &:= \{\kappa(|[\mathcal{F}f](t)|^2) \cdot \{m^{-1}|[\mathcal{F}f](t)|^{-2} \wedge 1\}\}, \qquad t \in \mathbb{R}. \end{aligned}$$

LEMMA A.4. *Suppose the samples from $f_Y$ and $f_\epsilon$ are independent. Let $f \in \mathcal{D}$, and consider $\mathcal{D}_f^\gamma$ defined in (4.2). Then, there exists $C > 0$, such that*

$$\inf_{\widehat{f_X}} \sup_{f_X \in \mathcal{M}_f^\rho, f_\epsilon \in \mathcal{D}_f^\gamma} \mathbb{E}\|\widehat{f_X} - f_X\|^2 \geq C \cdot \max_{t \in \mathbb{R}} \Delta_f^m(t).$$

PROOF. The proof is in analogy to the proof of Theorem 3.1 given in Neumann [31] and we omit the details. $\square$

PROOF OF THEOREM 4.2. The proof is based on the decomposition (4.1). From the bound given in (A.14), the polynomial source conditions (3.2) and (A.9) in Lemma A.2 with $w = \ell_s$ and $\gamma = \beta$, we obtain $\mathbb{E}\|\widehat{f_X}_s - \widehat{f_X^\alpha}_s\|_s^2 \leq \pi^{-1}\alpha^{-1} \cdot \mathbb{E}\|\widehat{f_Y} - f_Y\|^2 + 2C(\beta) \cdot \rho^2 \cdot \{\alpha^{-1} \cdot m^{-1-\beta} + (m \cdot \alpha)^{-1+\beta\wedge 1} \cdot m^{-\beta\wedge 1}\}$. While (A.15) in Lemma A.3 with $w = \ell_s$ and $\gamma = \beta$ provides $\mathbb{E}\|\widehat{f_X^\alpha}_s - f_X\|_s^2 \leq C_\beta \cdot \{\alpha^\beta + m^{-\beta}\} \cdot \rho^2$. Balancing these two terms then gives the result. $\square$

PROOF OF PROPOSITION 4.3. Let $g^s$ be defined by $\mathcal{F}g^s := \ell_s \cdot \mathcal{F}g$, $s \in \mathbb{R}$. Now, by making use of the relation $\|f_X^s\| = \|f_X\|_s$, the $H_s$-risk of an estimator $\widehat{f_X}$ of $f_X$ equals the MISE of $\widehat{f_X^s}$ as estimator of $f_X^s$. Moreover, $f_X$ belongs to $\mathcal{S}_f^\rho$ if and only if $f_X^s$ satisfies $\|\mathcal{F}f_X^s \cdot (|\mathcal{F}f^{-s}|^2)^{-\beta/2}\| \leq \rho$. Consider the sets $\mathcal{D}_f^\gamma$ and $\mathcal{M}_f^\rho$ defined in (4.2) and (A.24) with $k(t) = t^\beta$, respectively. Then, for any $f_0 \in \mathcal{D}_{f-s}^c$, $c > 0$, Lemma A.4 implies

$$\inf_{\widehat{f_X}} \sup_{f_X \in \mathcal{S}_f^\rho, f_\epsilon \in \mathcal{D}_f^\gamma} \mathbb{E}\|\widehat{f_X} - f_X\|_s^2 \geq \inf_{\widehat{f_X}} \sup_{f_X \in \mathcal{M}_{f_0}^{c\rho}, f_\epsilon \in \mathcal{D}_{f_0}^{c\gamma}} \mathbb{E}\|\widehat{f_X} - f_X\|^2$$

$$\geq C \max_{t \in \mathbb{R}} \left\{ |[\mathcal{F}f](t)|^{2\beta} \left\{ \frac{1}{m|[\mathcal{F}f](t)|^2} \wedge 1 \right\} \right\},$$

where the lower bound is of order $m^{-(1 \wedge \beta)}$, which proves the result. $\square$

PROOF OF COROLLARY 4.4. The proof is similar to the proof of Corollary 3.3, but uses Theorem 4.2 rather than Theorem 3.2, and we omit the details. $\square$

PROOF OF COROLLARY 4.5. Let $f \in \mathcal{D}_a$, and consider the sets $\mathcal{D}_f^\gamma$ and $\mathcal{S}_f^\rho$ defined in (4.2) and Proposition 4.3 with $\beta = (p-s)/(a+s)$, respectively. If $f_\epsilon \in \mathcal{D}_f^\gamma$, then $f_X \in H_p^\rho$ is equivalent to $f_X \in \mathcal{S}_f^{d\gamma\rho}$. Therefore, Proposition 4.3 leads to the following lower bound:

$$\inf_{\widehat{f_X}} \sup_{f_X \in H_p^\rho, f_\epsilon \in \mathcal{D}_a} \mathbb{E}\|\widehat{f_X} - f_X\|_s^2 \geq C m^{-(1 \wedge (p-s)/(a+s))}.$$

The result now follows by combination of the last lower bound with the lower bound in the case of known $f_\epsilon \in \mathcal{D}_a$ (cf. Mair and Ruymgaart [26]), that is, $\inf_{\widehat{f_X}} \sup_{f_X \in H_p^\rho, f_\epsilon \in \mathcal{D}_a} \mathbb{E}\|\widehat{f_X} - f_X\|_s^2 \geq C n^{-2(p-s)/(2(p+a)+1)}$. $\square$

PROOF OF THEOREM 4.6. Considering the decomposition (4.1), we bound the first term as in (A.12), and from (A.16) in Lemma A.3 with $w = \ell_s$ and $\gamma = \beta$, the second term satisfies $\mathbb{E}\|\widehat{f_{X_s}^\alpha} - f_X\|_s^2 \leq C_\beta |\log(C_\beta'\{\alpha + m^{-1}\})|^{-\beta}\rho^2$. The conditions of the theorem provide then $\mathbb{E}\|\widehat{f_{X_s}} - f_X\|_s^2 \leq C \cdot \{\mathbb{E}\|\widehat{f_Y} - f_Y\|^2 \vee m^{-1}\}^{1/2} + C \cdot |\log(C \cdot \{\mathbb{E}\|\widehat{f_Y} - f_Y\|^2 \vee m^{-1}\})|^{-\beta}$, for some constant $C > 0$ depending only on $\rho$ given in (3.2), $\beta$ and $c$, which implies the result. $\square$

PROOF OF PROPOSITION 4.7. The proof follows along the same lines as the proof of Proposition 4.3. Here, using the logarithmic rather than the polynomial source condition, Lemma A.4 implies

$$\inf_{\widehat{f_X}} \sup_{f_X \in \mathcal{S}_f^\rho, f_\epsilon \in \mathcal{D}_f^\gamma} \mathbb{E}\|\widehat{f_X} - f_X\|_s^2$$

$$\geq C \max_{t \in \mathbb{R}} \left\{ \frac{1}{|\log(|[\mathcal{F}f](t)|^2)|^\beta} \left\{ \frac{1}{m|[\mathcal{F}f](t)|^2} \wedge 1 \right\} \right\},$$

where the lower bound is of order $(\log m)^{-\beta}$, which gives the result. □

PROOF OF COROLLARY 4.8. The proof is similar to the proof of Corollary 3.5, but uses Theorem 4.6 rather than Theorem 3.4, and we omit the details. □

PROOF OF COROLLARY 4.9. The proof follows along the same lines as the proof of Corollary 4.5. Here, using Proposition 4.7 rather than Proposition 4.3 leads to the lower bound $C(\log m)^{-(p-s)/a}$. The result follows then from the lower bound $C(\log n)^{-(p-s)/a}$ in the case of known $f_\epsilon$ (cf. Mair and Ruymgaart [26]). □

PROOF OF THEOREM 4.10. The proof is based on the decomposition (4.1). From the bound given in (A.14), the general source conditions (3.6) and (A.10) in Lemma A.2 with $w = \ell_s$, we obtain $\mathbb{E}\|\widehat{f_X}_s - \widehat{f_X^\alpha}_s\|_s^2 \leq \pi^{-1}\alpha^{-1}\mathbb{E}\|\widehat{f_Y} - f_Y\|^2 + 2C\rho^2\{\alpha^{-1}m^{-1}\kappa(m^{-1}) + \kappa(m^{-1})\}$. While (A.17) in Lemma A.3 with $w = \ell_s$ provides $\mathbb{E}\|\widehat{f_X^\alpha}_s - f_X\|_s^2 \leq C_\kappa\kappa(C_\kappa\{\alpha + m^{-1}\})\rho^2$. The condition on $\alpha$ ensures then the balance of these two terms. The result follows by making use of the relation $\kappa(c \cdot \delta) \leq (c \vee 1) \cdot \kappa(\delta)$, which follows, for $c < 1$ and for $c \geq 1$, from the monotonicity and the concavity of $\kappa$, respectively. □

PROOF OF PROPOSITION 4.11. The proof follows along the same lines as the proof of Proposition 4.3. Here, using the general rather than the polynomial source condition, Lemma A.4 implies

$$\inf_{\widehat{f_X}} \sup_{f_X \in \mathcal{S}_f^\rho, f_\epsilon \in \mathcal{D}_f^\gamma} \mathbb{E}\|\widehat{f_X} - f_X\|_s^2 \geq C \max_{t \in \mathbb{R}}\left\{\kappa(|[\mathcal{F}f](t)|^2)\left\{\frac{1}{m|[\mathcal{F}f](t)|^2} \wedge 1\right\}\right\}.$$

Since $\kappa$ is increasing and $\kappa(t^2)/t^2$ is decreasing, it follows that the lower bound is of order $\kappa(1/m)$, which proves the result. □

**Acknowledgments.** I thank the referees for their careful reading of the paper and for helpful comments and suggestions.

## REFERENCES

[1] BIGOT, J. and VAN BELLEGEM, S. (2006). Log-density deconvolution by wavelet thresholding. Technical report, Univ. catholique de Louvain.
[2] BUTUCEA, C. and MATIAS, C. (2005). Minimax estimation of the noise level and of the deconvolution density in a semiparametric deconvolution model. *Bernoulli* **11** 309–340. MR2132729
[3] CARRASCO, M. and FLORENS, J.-P. (2002). Spectral method for deconvolving a density. Working Paper No. 138, Univ. Toulouse I, IDEI.
[4] CARROLL, R. J. and HALL, P. (1988). Optimal rates of convergence for deconvolving a density. *J. Amer. Statist. Assoc.* **83** 1184–1186. MR0997599

[5] CAVALIER, L. and HENGARTNER, N. (2005). Adaptive estimation for inverse problems with noisy operators. *Inverse Problems* **21** 1345–1361. MR2158113

[6] COMTE, F., ROZENHOLC, Y. and TAUPIN, M.-L. (2006). Penalized contrast estimator for density deconvolution. *Canad. J. Statist.* **34** 431–452. MR2328553

[7] DEVROYE, L. (1989). Consistent deconvolution in density estimation. *Canad. J. Statist.* **17** 235–239. MR1033106

[8] DIGGLE, P. J. and HALL, P. (1993). A Fourier approach to nonparametric deconvolution of a density estimate. *J. Roy. Statist. Soc. Ser. B* **55** 523–531. MR1224414

[9] EFROMOVICH, S. (1997). Density estimation for the case of supersmooth measurement error. *J. Amer. Statist. Assoc.* **92** 526–535. MR1467846

[10] ENGL, H. W., HANKE, M. and NEUBAUER, A. (2000). *Regularization of Inverse Problems.* Kluwer Academic, Dordrecht. MR1408680

[11] FAN, J. (1991). On the optimal rates of convergence for nonparametric deconvolution problems. *Ann. Statist.* **19** 1257–1272. MR1126324

[12] FAN, J. (1992). Deconvolution with supersmooth distributions. *Canad. J. Statist.* **20** 734–747. MR1183078

[13] FAN, J. and KOO, J. Y. (2002). Wavelet deconvolution. *IEEE Trans. Inform. Theory* **48** 734–747. MR1889978

[14] GOLDENSHLUGER, A. (1999). On pointwise adaptive nonparametric deconvolution. *Bernoulli* **5** 907–925. MR1715444

[15] GOLDENSHLUGER, A. (2000). Density deconvolution in the circular structural model. *J. Multivariate Anal.* **81** 360–375. MR1906385

[16] HALL, P. and QIU, P. (2005). Discrete-transform approach to deconvolution problems. *Biometrika* **92** 135–148. MR2158615

[17] HALL, P. and YAO, Q. (2003). Inference in components of variance models with low replication. *Ann. Statist.* **31** 414–441. MR1983536

[18] HAVILIO, M. (2006). Signal deconvolution based expression-detection and background adjustment for microarray data. *J. Comput. Biol.* **13** 63–80. MR2253541

[19] HOROWITZ, J. L. and MARKATOU, M. (1996). Semiparametric estimation of regression models for panel data. *Rev. Econom. Stud.* **63** 145–168. MR1372250

[20] KAWATA, T. (1972). *Fourier Analysis in Probability Theory.* Academic Press, New York. MR0464353

[21] KIM, P. T. and KOO, J. Y. (2002). Optimal spherical deconvolution. *J. Multivariate Anal.* **80** 21–42. MR1889831

[22] KOO, J. Y. and PARK, B. U. (1996). $B$-spline deconvolution based on the EM algorithm. *J. Statist. Comput. Simulation* **54** 275–288. MR1701219

[23] LEVITT, D. G. (2003). The use of a physiologically based pharmacokinetic model to evaluate deconvolution measurements of systemic absorption. In *BMC Clinical Pharmacology* **3**. BioMed Central. Available at http://www.biomedcentral.com/1472-6904/3/1.

[24] LIU, M. C. and TAYLOR, R. L. (1989). A consistent nonparametric density estimator for the deconvolution problem. *Canad. J. Statist.* **17** 427–438. MR1047309

[25] MAIR, B. A. (1994). Tikhonov regularization for finitely and infinitely smoothing operators. *SIAM J. Math. Anal.* **25** 135–147. MR1257145

[26] MAIR, B. A. and RUYMGAART, F. H. (1996). Statistical inverse estimation in Hilbert scales. *SIAM J. Math. Anal.* **56** 1424–1444. MR1409127

[27] MEISTER, A. (2006). Density estimation with normal measurement error with unknown variance. *Statist. Sinica* **16** 195–211. MR2256087

[28] MENDELSOHN, J. and RICE, J. (1982). Deconvolution of micro-fluorometric histograms with $B$-splines. *J. Amer. Statist. Assoc.* **77** 748–753.

[29] NAIR, M., PEREVERZEV, S. V. and TAUTENHAHN, U. (2005). Regularization in Hilbert scales under general smoothing conditions. *Inverse Problems* **21** 1851–1869. MR2183654
[30] NATTERER, F. (1984). Error bounds for Tikhonov regularization in Hilbert scales. *Applicable Anal.* **18** 29–37. MR0762862
[31] NEUMANN, M. H. (1997). On the effect of estimating the error density in nonparametric deconvolution. *J. Nonparametr. Statist.* **7** 307–330. MR1460203
[32] NEUMANN, M. H. (2006). Deconvolution from panel data with unknown error distribution. *J. Multivariate Anal.* **98** 1955–1968. MR2396948
[33] PARZEN, E. (1962). On estimation of a probability density function and mode. *Ann. Math. Statist.* **33** 1065–1076. MR0143282
[34] PENSKY, M. and VIDAKOVIC, B. (1999). Adaptive wavelet estimator for nonparametric density deconvolution. *Ann. Statist.* **27** 2033–2053. MR1765627
[35] PETROV, V. V. (1995). *Limit Theorems of Probability Theory. Sequences of Independent Random Variables*, 4th ed. Clarendon Press, Oxford. MR1353441
[36] STEFANSKI, L. A. (1990). Rate of convergence of some estimators in a class of deconvolution problems. *Statist. Probab. Lett.* **9** 229–235. MR1045189
[37] TAUTENHAHN, U. (1996). Error estimates for regularization methods in Hilbert scales. *SIAM J. Numer. Anal.* **33** 2120–2130. MR1427456
[38] TESSIER, E. (1995). Analysis and calibration of natural guide star adaptive optics data. In *Adaptive Optical Systems and Applications* (R. K. Tyson and R. Q. Fugate, eds.) 178–193. *Proc. SPIE* **2534**. Available at arXiv:astro-ph/9601174v1.
[39] VAINIKKO, G. M. and VERETENNIKOV, A. Y. (1986). *Iteration Procedures in Ill-Posed Problems*. Nauka, Moscow (in Russian). MR0859375
[40] VAN DER VAART, A. W. (1998). *Asymptotic Statistics.* Cambridge Univ. Press. MR1652247
[41] ZHANG, C.-H. (1990). Fourier methods for estimating mixing densities and distributions. *Ann. Statist.* **18** 806–831. MR1056338

INSTITUTE OF APPLIED MATHEMATICS
RUPRECHT–KARLS–UNIVERSITÄT HEIDELBERG
IM NEUENHEIMER FELD 294
D-69120 HEIDELBERG
GERMANY
E-MAIL: johannes@statlab.uni-heidelberg.de